\documentclass{amsart}

\usepackage[all]{xy}

\newfont{\fnt}{bbold10 scaled 1250}

\newcommand{\J}{^{\scriptscriptstyle J}\displaystyle\!}
\newcommand{\Oo}{\mathcal O}
\newcommand{\Dd}{\mbox{$\not\!\partial$}}

\newcommand{\Gg}{\mathcal G}
\newcommand{\Ee}{\mathcal E}

\newcommand{\Mm}{\mathfrak M}

\newcommand{\Cc}{\mathbb C}

\newcommand{\Zz}{\mathbb Z} 

\newcommand{\Id}{\mbox{\fnt 1}}

\newcommand{\CP}{\mathbb{CP}}

\newtheorem{teor}{Theorem}[section]
\newtheorem{defi}[teor]{Definition}

\newtheorem{prop}[teor]{Proposition}
\newtheorem{conj}[teor]{Conjecture}

\title[Dirac operator coupled to instantons]{Dirac operator coupled to
Instantons on positive definite 4 manifolds}

\author{Jo\~ao Paulo Santos}

\begin{document}

\bibliographystyle{amsalpha}

\setlength{\parindent}{0em}
\setlength{\parskip}{1.5ex}
\newdir{ >}{{}*!/-5pt/\dir{>}}
\newdir{ |}{{}*!/-5pt/\dir{|}}

\begin{abstract}
We study the moduli space of instantons on a simply connected positive 
definite four manifold 
by analyzing the classifying map of the index bundle of a family of Dirac
operators parametrized by the moduli space. As applications we compute the
cohomology ring for the charge 2 moduli space in the rank stable limit.
\end{abstract}

\maketitle

\section{Introduction}

Let $X$ be a simply connected positive definite four manifold 
with base point $x_\infty\in X$ and let
$E\to X$ be an $SU(r)$ vector bundle with second Chern class $c_2(E)=-k<0$.
We will denote by $\Mm_k^r(X)$ the moduli space of self-dual
connections on $E$, framed at $x_\infty$. By results in \cite{Don83} and
\cite{Fre82}, $X$ is homeomorphic to
a connected sum of a number $q$ of copies of $\CP^2$. Let $X_q$ denote this
manifold.

A Dirac operator $\Dd:\Gamma(S^+)\to\Gamma(S^-)$ defines
a $K$-theory class $[Ker\,\Dd]-[Coker\,\Dd]$ over $\Mm_k^r(X_q)$
whose fiber at each point $A\in\Mm_k^r(X_q)$
is given by the kernel and cokernel of
the operator $\Dd$ coupled to the connection $A$.
For $A$ a self-dual connection, the operator $\Dd$ is surjective,
hence its kernel defines a vector bundle over $\Mm_k^r(X_q)$,
with classifying map $f:\Mm_k^r(X_q)\to BU(N)$
($N={\rm dim\,Ker}\,\Dd={\rm Ind}\,\Dd$).
In this paper we will study this family of maps 
(one for each Dirac operator) and use it to
obtain information about the topology of the moduli space $\Mm_k^r(X_q)$.

For $r'>r$, direct sum with a trivial rank $r'-r$ bundle induces a map
$\Mm_k^r(X_q)\to\Mm_k^{r'}(X_q)$ and we define 
$\displaystyle\Mm_k^\infty(X_q)=\lim_{\substack{\longrightarrow\\r}}
\Mm_k^r(X_q)$.
In \cite{San95} it was shown that, for $q=0$
($X_0=S^4$), ${\rm Ind}\,\Dd=k$ and the map
$f:\Mm_k^\infty(S^4)\to BU(k)$ 
is a homotopy equivalence. Our first result is
(see also \cite{BrSa97})

\begin{teor}\label{t1}
Let $\Dd_\pm$ be the Dirac operators on $\CP^2$
with associated Chern classes the generators $\pm1\in H^2(\CP^2)$. Then 
${\rm Ind}\,\Dd_\pm=k$ and the associated map
$f_+\times f_-:\Mm_k^\infty(\CP^2)\to BU(k)\times BU(k)$
is a homotopy equivalence.
\end{teor}

These results have a nice interpretation as follows:
In \cite{Tau89} Taubes showed that, for $k'>k$, there is
a glueing map $\Mm_k^r(X)\to\Mm_{k'}^r(X)$ such that the space
$\displaystyle\Mm_{\infty}^r(X)=\lim_{\substack{\longrightarrow\\k}}
\Mm_k^r(X)$
is homotopy equivalent to ${\rm Map}(X,BSU(r))$. When $r\to\infty$
${\rm Map}(X_q,BSU)$ is homotopic to the product of $q+1$ copies of $BU$
so we have a map 
$\Mm_k^\infty(X_q)\to\displaystyle
\lim_{\substack{\longrightarrow\\k}}\Mm_k^\infty\cong\prod_{l=0}^q BU$.
For $X_0,X_1$ this map is just the map $BU(k)\to BU$ and
$BU(k)\times BU(k)\to BU\times BU$ respectively. For each $q$
we will choose $q+1$
Dirac operators with index $k$
inducing a map $f:\Mm_k^r(X_q)\to\displaystyle
\prod_{l=0}^q BU(k)$ such that the diagram
\[
\xymatrix{
\Mm_k^r(X_q)\ar[r]^-{f}\ar[d]&\prod_{l=0}^q BU(k)\ar[d]\\
\Mm_\infty^r(X_q)\ar[r]^-{f}&\prod_{l=0}^q BU
}\]
commutes and, when $r\to\infty$, the bottom map is
a homotopy equivalence.
We will show that

\begin{teor}\label{t3}
For $k=1,2$ and any $q$, the map 
$\Mm_k^\infty(X_q)\to\displaystyle
\lim_{\substack{\longrightarrow\\k}}\Mm_k^\infty\cong\textstyle\prod BU$
induces a surjective map in cohomology. 
\end{teor}

We conjecture this result is true for all $k$:

\begin{conj}
The map
$\Mm_k^\infty(X_q)\to\displaystyle
\lim_{\substack{\longrightarrow\\k}}\Mm_k^\infty\cong\textstyle\prod BU$
induces a surjective map in cohomology.
\end{conj}

Theorem \ref{t3}
will allow us to compute the cohomology ring of $\Mm_2^\infty(X_q)$:

\begin{teor}\label{t2}
$H^*(\Mm_2(X_q))$ is isomorphic to the quotient of the polynomial ring
$\Zz[C_1,C_2,S_1^i,S_2^i;\,1\leq i\leq q]$ 
(with deg $C_n,S_n^i=2n$) by the ideal generated by
\begin{align*}
&C_1S_1^iS_1^j+S_1^iS_2^j+S_2^iS_1^j\,,\ i\neq j\\
&C_2S_1^iS_1^j-S_2^iS_2^j\,,\ i\neq j\\
&S_1^iS_1^jS_1^k\,,\ i\neq j\neq k\\
&S_1^iS_1^jS_2^k\,,\ i\neq j\neq k
\end{align*}
\end{teor}

The strategy of the proofs will be
to define a subspace $\Sigma_{q,k}^r\subset\Mm_k^r(X_q)$ and 
to compute the 
map induced in cohomology by the composition
\[
\Sigma_{q,k}^r\to\Mm_k^r(X_q)\to\prod_{l=0}^q BU(k)
\]
In section 2 we will define the moduli spaces $\Mm_k^r(X_q)$ and
the classifying map
$f:\Mm_k^r(X_q)\to\prod BU(k)$. Then in section 3 we will introduce the
spaces $\Sigma_{q,k}^r$ and use them to prove theorem \ref{t1}.
In section 4 we define, for $q'<q$, maps $\pi^*:\Mm_k^r(X_{q'})\to
\Mm_k^r(X_q)$ and study some of their properties. These results will
be used in section 5 to prove theorems \ref{t2} and \ref{t3}.

\section{The moduli space and the classifying map}

Let $X_q$ denote the connected sum of $q$ copies of $\CP^2$.
Fix a $C^\infty$ $SU(r)$ vector bundle $E\to X_q$ with $c_2(E)=-k<0$
and let $\mathcal A(E)$ denote the space of connections on $E$.
A connection $A\in\mathcal A(E)$ is called self-dual if
$F_A^-=0$. These connections minimize the Yang-Mills functional
$\int |F_A|^2$ (see \cite{DoKr90}).
Let $\mathcal G$ be the gauge group of automorphisms of the bundle $E$.
Fix a point $x_\infty\in X_q$ and let
$\mathcal G_0=\left\{g\in\Gg\,|g(x_\infty)=\Id\right\}$.
The moduli space $\Mm_k^r(X_q)$
of instantons framed at $x_\infty$ is the quotient of the space
of self-dual connections by $\mathcal G_0$.

Given a $Spin^c$ structure on $X_q$ with associated line bundle $L$
and cohomology class $c\in H^2(X_q)$
($c=w_2(X_q)\,{\rm mod}\,2$),
a choice of a connection $a$ in $L$ induces a Dirac operator
$\Dd_c:\Gamma(S^+)\to\Gamma(S^-)$ (see for example \cite{LaMi89}). 
We will take $a$ to be the (unique) self-dual connection on $L$.
Then, to each connection 
$A\in\mathcal A(E)$ we can associate a Dirac operator 
$\Dd_{c,A}:\Gamma(S^+\otimes E)\to\Gamma(S^-\otimes E)$
coupled with $A$. 

\begin{prop}
The index of $\Dd_{c,A}$ is given by
\[
ind\,\Dd_{c,A}=k+\frac{c_1(E)\cdot(c_1(E)+c)}2+r\frac{c^2-q}8
\]
\end{prop}
\begin{proof}
It is a direct consequence of the formula
\[
ind\,\Dd_{c,A}=e^{\frac c2}\cdot ch(E)\cdot\hat A(M)
\]
for the index of a Dirac operator.
\end{proof}

Notice that if $c^2={\rm sig}(M)=q$
then $ind\,\Dd_{c,A}$ is independent of the rank $r$ of $E$. 

\begin{prop}
There are
exactly $2^q$ $Spin^c$ structures on $M$ such that
$c^2=q$, namely $c=(\pm1,\ldots,\pm1)$.
\end{prop}
\begin{proof}
From Wu's formulas,
\[
w_2=Sq^0v_2+Sq^1v_1+Sq^2v_0=v_2
\]
where $v_2\in H^2(M,\Zz_{/2})$ satisfies $v_2\cup x=Sq^2x=x\cup x$
for any $x\in H^2(M,\Zz_{/2})$ (see \cite{MiSt74}, theorem 11.14). 
From here it follows that
\[
w_2(M)=(1,\ldots,1)\in H^2(M,\Zz_{/2})
\]
Now the conditions $c^2={\rm sig}(M)=q$ and $w_2=c\,{\rm mod}\,2$ imply
$c=(\pm1,\ldots,\pm1)$.
\end{proof}

\begin{prop}
Let $X_q,\Dd_{c,A}$ be as above and assume that $X_q$ has a metric with
scalar curvature $s>0$.
Then if the connections $A,a$ are
self-dual, the operator $\Dd_{c,A}$ is surjective.
\end{prop}
\begin{proof}
We want to show the cokernel of $\Dd_{c,A}$ is zero. Hence we look at
the Bochner-Weizenboch formula for the dual operator
$\Dd_{c,A}^*$:
\[
\Dd_{c,A}\Dd_{c,A}^*=\nabla^*\nabla+\frac14s+\mathfrak R
\]
where $s$ is the scalar curvature of $X_q$ and
$\mathfrak R$ is defined as follows: locally we can write
$S^-=S_0\otimes L^{\frac12}\otimes E$. Let $\sigma\in\Gamma(S_0)$,
$v\in\Gamma(L^{\frac12}\otimes E)$. Then
\[
\mathfrak R(\sigma\otimes v)=\sum_{j,k=1}^4
(e_je_k\sigma)\otimes(F^-(e_j,e_k)v)
\]
Here $F^-$ is the anti-self-dual part of the curvature on 
$L^{\frac12}\otimes E$ induced by $a,A$, which vanishes if 
$A,a$ are self-dual. Since $s>0$ it follows that $\Dd_{c,A}^*$ is injective
which completes the proof.
\end{proof}

Let $\mathcal I\subset\mathcal A$ denote the space of self-dual connections,
We have the diagram
\[
\xymatrix{
\mathcal I\times_{\Gg}\Gamma(S^+\otimes E)\ar[rr]^-{\Dd_c}\ar[rd]_{p_1}&&
\mathcal I\times_{\Gg}\Gamma(S^-\otimes E)\ar[ld]^{p_2}\\
&\mathcal I/\mathcal G_0}
\]
where $\Dd_c([A,s])=[A,\Dd_{c,A}s]$. Then the kernel of
$\Dd_c$ defines a rank $k$ vector bundle $[{\rm Ker}\,\Dd_c]$ over 
$\mathcal I/\mathcal G_0=\Mm_k^r$ (see \cite{Seg70}),
with classifying map $f_c:\Mm_k^r\to BU(k)$.

Now we choose $q+1$ 
specific $Spin^c$ structures: fix a basis $e_1,\ldots,e_q$ of $H^2(X_q)$.
Then define
$c_0=(1,1,\ldots,1)$ and, for $l=1,\ldots,q$,
\[
c_l=-e_l+\sum_{i\neq l}e_i
\]
The associated classifying maps define a map
$f:\Mm_k^r(X_q)\to\displaystyle\prod_{l=0}^qBU(k)$.
This choice of Chern classes $c_0,\ldots,c_q$
is justified by the following theorem
(which will be proven at the end of next section):

\begin{teor}\label{t4}
When $k,r\to\infty$, the map $f:\Mm_\infty^\infty(X_q)\to\prod BU$
is a homotopy equivalence.
\end{teor}

\section{Orbits}

In this section we define the subspace 
$\Sigma_{q,k}^r\subset\Mm_k^r(X_q)$
and compute the composition 
$\Sigma_{q,k}^r\to\Mm_k^r\to\prod BU(k)$.
As a corollary we prove theorems \ref{t1} and \ref{t4}.

The group $\mathcal G$ of automorphisms of $E$
acts on $\Mm_k^r(X_q)$.
For a connection $A\in\mathcal A(E)$ let 
$Or(A)\subset\Mm_k$ denote the orbit of $[A]$ under the action of $\Gg$. 
We will choose $\binom{k+q-1}{q-1}$
connections $A_J$ and define $\Sigma_{q,k}^r$ as the
union of the orbits $Or(A_J)$:

To each q-tuple $J=(j_1,j_2,\ldots,j_q)$ with
$|J|=\sum_ij_i=k$ we define a connection $A_J$ on $E$ as follows:
fix an isomorphism
$\phi:E\cong\bigoplus_aL_a$ where the line bundles $L_a$ satisfy
\[
\#\{a|c_1(L_a)=e_i\}=\#\{a|c_1(L_a)=-e_i\}=j_i\ ,\
\#\{a|c_1(L_a)=0\}=r-2k
\]
($e_1,\ldots,e_q$ is a basis of $H^2$).
That is, $E\cong L_{+e_1}^{\oplus j_1}\oplus
L_{-e_1}^{\oplus j_1}\oplus\ldots\oplus
L_{+e_q}^{\oplus j_q}\oplus L_{-e_q}^{\oplus j_q}\oplus L_0^{r-2k}$.
Let $\lambda_a$ be the self-dual connection on $L_a$ and define
$A_J=\phi^{-1}\left(\bigoplus_a\lambda_a\right)\phi$.
We can now define $\Sigma_{q,k}^r$:

\begin{defi}
The space $\Sigma_{q,k}^r$ is the disjoint union 
\[
\Sigma_{q,k}^r=\coprod_{|J|=k} Or(A_J).
\]
\end{defi}

Observe that
for a different choice of isomorphism $\hat\phi:E\cong\bigoplus_aL_a$, 
$A_J$ would differ by $\hat\phi^{-1}\phi\in\Gg$.
Hence $Or(A_J)$ depends only on $J$.
We will often simplify the notation and write $Or(J)$
instead of $Or(A_J)$.

\begin{prop}\label{prop31}
Let $I_A$ be the stabilizer of $A$ under the action of $\Gg$,
$I_{[A]}$ be the stabilizer of $[A]$ under $\Gg/\Gg_0$ and
$[I_A]$ the image of $I_A$ on $\Mm_k$. Then these 3 groups
are isomorphic and we have
\[
Or(A)\cong\frac{SU(r)}{I_A}
\]
\end{prop}
\begin{proof}
It follows directly
from the fact that the action of $\Gg_0$ on $\mathcal A(E)$ is free.
\end{proof}

\begin{prop}
The stabilizer $I_{[A_J]}\subset\Gg/\Gg_0$
of $[A_J]$ is the intersection of the
subgroup
\[
P(J)=U(j_1)\times U(j_1)\times
\ldots\times U(j_q)\times U(j_q)\times U(r-2k),
\]
sitting in $U(r)$ diagonally, with $SU(r)$.
\end{prop}
\begin{proof}
It is clear that $P(J)\subset I_{[A_J]}$.
Let $g\in I_{A_J}\subset\Gg$ (we are identifying $I_{A_J}$ and
$I_{[A_J]}$: see proposition \ref{prop31}).
Locally $(A_J)_{ab}=\delta_{ab}\lambda_a$ and we get the equation
$dg_{ab}=g_{ab}(\lambda_a-\lambda_b)$ which implies
$g_{ab}(d\lambda_a-d\lambda_b)=0$. So, since $d\lambda_a$ is the
harmonic 2-form representing $c_1(L_a)$,
$g_{ab}=0$ if $c_1(L_a)\neq c_1(L_b)$. This shows $g\in P(J)$ completing
the proof.
\end{proof}

Hence, when $r\to\infty$, the orbit of $A_J$ is given by
\[
Or(J)\cong BU(j_1)^{\times2}\times\ldots\times BU(j_q)^{\times2}
\]
The main result of this section is

\begin{teor}\label{teor34}
Let
\[
B\pi_c:BU(j_1)^{\times2}\times\ldots\times BU(j_q)^{\times2}\to
BU(j_1)\times\ldots\times BU(j_q)
\]
be the projection map defined as follows: in each factor $BU(j_i)^{\times2}$,
$B\pi_c$ is the projection onto the first component if $c\cup e_i=1$
and onto the second component if $c\cup e_i=-1$. Then the restriction of
$f$ to $Or(J)$, $f:Or(J)\to BU(k)$,
is the composition of $B\pi_c$ with the Whitney sum map
\[
BU(j_1)\times\ldots\times BU(j_q)\to BU(k)
\]
\end{teor}

Before we prove this theorem we need the following lemma
analyzing the restriction of the vector space
$[{\rm Ker}\,\Dd_c]\to\Mm_k^r(X_q)$ to an orbit $Or(A_0)$:

\begin{prop}
Let $G=\Gg/\Gg_0$ and fix $A_0\in\mathcal I$.
Let $\rho:I_{[A_0]}\to I_{A_0}$ denote the isomorphism of
proposition \ref{prop31}.
Then $\rho$ induces a representation of $I_{[A_0]}$ on the vector space
\[
V=\{\psi\in\Gamma(S^+\otimes E)\,|\,\Dd_{A_0}\psi=0\}
\]
The bundle $[Ker\,\Dd_c]$ restricted to the orbit $Or(A_0)$
is isomorphic to the vector bundle associated to the principal bundle
$G\to G/I_{[A_0]}$ and the representation $\rho$.
\end{prop}
\begin{proof}
$\rho$ is well defined since, for $h\in I_{A_0}$,
$\Dd_{A_0}\psi=h^{-1}\Dd_{A_0}h\psi$.
We begin by defining a map
\[
\xymatrix@R=0em{
\Gg\times V\ar[r]&\mathcal I\times_{\Gg_0}\Gamma\\
(g,\psi)\ar@{|->}[r]&[g^{-1}A_0g,g^{-1}\psi]
}\]
(where $g^{-1}A_0g$ is a shorthand for 
$g^{-1}\circ\nabla_{A_0}\circ g$).
Observe that this map descends to give a map
$(\Gg/\Gg_0)\times_{I_{A_0}} V\to\mathcal I\times_{\Gg_0}\Gamma$ where
$\Gg_0$ acts on $\Gg$ on the right and $I_{A_0}$ acts on $\Gg$ on the left
and on $V$ by $\rho$.
The statement of the proposition is that this bundle map,
when restricted to $Or(A)$, is an isomorphism onto
the kernel of $\Dd_c$. 

We first show surjectivity. Let $A\in Or(A),\,
\hat\psi\in\Gamma(S^+\otimes E)$ be such that
$\Dd_{A}\hat\psi=0$. Then $A=g^{-1}\circ A_0\circ g$ for some
$g\in\Gg$ and $g\hat\psi\in V$ because
$\Dd_{A_0}g\hat\psi=g\Dd_A\psi=0$. Then
$(g,g\hat\psi)\mapsto[A,\hat\psi]$ which shows surjectivity.

Now suppose $(g_1,\psi_1)$ and $(g_2,\psi_2)$ have the same image, 
that is, 
$[g_1^{-1}A_0g_1,g_1^{-1}\psi_1]=
[g_2^{-1}A_0g_2,g_2^{-1}\psi_2]$. Then, for some
$h\in \Gg_0$, $(g_1h)^{-1}\psi_1=g_2^{-1}\psi_2$
and $(g_1h)^{-1}A_0(g_1h)=g_2^{-1}A_0g_2$, so
$g_2(g_1h)^{-1}\in I_{A_0}$. It follows that
\[
(g_1,\psi_1)\sim
(g_1h,\psi_1)\sim
(g_2(g_1h)^{-1}g_1h,g_2(g_1h)^{-1}\psi_1)=(g_2,\psi_2)
\]
This concludes the proof of the proposition.
\end{proof}

We can now prove theorem \ref{teor34}:

\begin{proof}
$f$ is the composition of the classifying map
$G/I_{A_J}\to BI_{A_J}$ with the map
$B\rho:BI_{A_J}\to BGl(V)$ induced by the
representation $\rho:I_{A_J}\to Gl(V)$.
So we have to analyze this representation.
Let $\psi\in V$ and let $g\in P(J)$. Recall
\[
E\cong L_{+e_1}^{\oplus j_1}\oplus L_{-e_1}^{\oplus j_1}\oplus\ldots\oplus
L_{+e_q}^{\oplus j_q}\oplus L_{-e_q}^{\oplus j_q}\oplus L_0^{\oplus r-2k}
\]
Then we can write
\begin{align*}
&\psi=\psi_1^+\oplus\psi_1^-\oplus\ldots\oplus\psi_q^+\oplus\psi_q^-\oplus\psi_0\\
&g   =   g_1^+\oplus   g_1^-\oplus\ldots\oplus   g_q^+\oplus   g_q^-\oplus g_0
\end{align*}
with $\psi_i^\pm\in\Gamma\left(L_{\pm e_i}^{\oplus j_i}\right)$
and equivalently for $g_i^\pm$.
Now $\Dd_A\psi=0$ hence $\Dd_{A_i^\pm}\psi_i^\pm=0$ where
$A_i^\pm=\lambda_{\pm e_i}^{\oplus j_i}$
is the connection on $L_{\pm e_i}^{\oplus j_i}$. 
Also $\Dd\psi_0=0$ hence $\psi_0=0$ and $g_0$ acts trivially. From the index
theorem it follows
\[
ind\,\Dd_{\lambda_{\pm e_i}}=\frac{\pm e_i\cdot(c\pm e_i)}2=
\frac{1\pm c\cup e_i}2
\]
Now, $c\cup e_i=\pm1$. Hence
\begin{itemize}
\item If $c\cup e_i=+1$ then $\psi_i^-=0$ and
$\psi_i^+=v\psi_i$ with $\psi_i\in\Gamma(L_{+e_i})$ and $v\in\Cc^{j_i}$;
It follows that $g_i^-$ acts trivially and $g_i^+$ acts by matrix
multiplication on $v$.
\item If $c\cup e_i=-1$ then $\psi_i^+=0$ and
$\psi_i^-=v\psi_i$ with $\psi_i\in\Gamma(L_{-e_i})$ and $v\in\Cc^{j_i}$;
It follows that $g_i^+$ acts trivially and $g_i^-$ acts by matrix
multiplication on $v$.
\end{itemize}
Fixing $\psi_1,\ldots,\psi_q$ as above gives an isomorphism
$V\cong\Cc^{j_1}\oplus\ldots\oplus\Cc^{j_q}$.
Let
\[
P_c(J)\cong\left\{U(j_1)\times\ldots\times U(j_q)\,|\,{\rm det}=1\right\}
\]
Then there is an obvious representation $\rho_c$ of $P_c(n)$ on
$V\cong \Cc^{j_1}\oplus\ldots\oplus\Cc^{j_q}$ which is trivial on the
$U(r-k)$ factor. It follows that $\rho$ is the composition of $\pi_c$
with direct sum $P_c\to Gl(V)$.

Now, when $r\to\infty$ we have
\[
G/I_{A_J}\cong
BU(j_1)^{\times2}\times\ldots\times BU(j_q)^{\times2}
\]
and $BI_{A_J}\cong G/I_{A_J}\times BU^{\times 2}$.
The classifying map of $G\to G/I_{A_J}$ is
then the identity on the $G/I_{A_J}$ component
and the result follows.
\end{proof}

As a corollary we have

\begin{teor}
For $X_1=\CP^2$, and when $r\to\infty$,
the map $f:\Mm_k^\infty(X_1)\to BU(k)\times BU(k)$ 
is a homotopy equivalence.
\end{teor}
\begin{proof}
Clearly the restriction of $f$ to $\Sigma_{1,k}^\infty$ is a
homotopy equivalence. From theorem \ref{teorA2}
in the appendix,
the inclusion $\Sigma_{1,k}^\infty\to\Mm_k^\infty(X_1)$
is a homotopy equivalency.
The result follows.
\end{proof}

As another corollary we can now prove theorem \ref{t4}:

\begin{teor}
When $k,r\to\infty$, the map $f:\Mm_\infty^\infty(X_q)\to\prod BU$
is a homotopy equivalence.
\end{teor}
\begin{proof}
We begin by writing $\Mm_\infty^\infty=\displaystyle
\lim_{\substack{\longrightarrow\\n}}\Mm_{nq}^\infty$.
Let $[A_n]\in\Mm_{nq}^\infty$ be the connection 
associated with $J=(n,n,\ldots,n)$.
We claim that, when $n\to\infty$, the map
\[
f_*:H_*\left(\,Or(A_n)\,\right)\to
H_*\left(\,BU(nq)\times\ldots\times BU(nq)\,\right)
\]
induced by $f$ in homology is surjective.
It follows that, when $n\to\infty$, the map 
\[
f_*:H_*\left(\,\Mm_{nq}^\infty\,\right)\to
H_*\left(\,BU(nq)\times\ldots\times BU(nq)\,\right)
\]
is surjective hence an isomorphism. Since $\Mm_\infty^\infty$ is
simply connected,
$f$ is a homotopy equivalence.

To prove the claim we recall that
\[
Or(A_n)\cong \prod_{i=1}^q\left(\, BU(n)\times BU(n)\,\right)
\]
When $n\to\infty$, $Or(A_\infty)$ is a loop space and 
the homology ring is
\begin{align*}
&H_*(Or(A_\infty))\cong \Zz[x_m^{i+},x_m^{i-};1\leq i\leq q,m\geq1]\\
&H_*\left(\prod_{l=0}^q BU\right)\cong \Zz[X_m^l;0\leq l\leq q,m\geq 1]
\end{align*}
and, by theorem \ref{teor34}, $f_*$ is given by
\begin{align*}
&f_*(x_m^{i+})=\left(\sum_{l=0}^qX_m^l\right)-X_m^i\\
&f_*(x_m^{i-})=X_m^i
\end{align*}
This map is clearly surjective. This completes the proof.
\end{proof}

\section{Instantons and holomorphic bundles}

The objective of this section is to define, for $q'<q$,
embeddings $\pi^*:\Mm_k^r(X_{q'})\to\Mm_k^r(X_q)$ and
to investigate their properties. In particular we will
show that 
\begin{itemize}
\item $\pi^*(Or(A))=Or(\pi^*A)$
\item For $X_{q'}=S^4$ and $r\to\infty$, $f_c\circ\pi^*$
is a homotopy equivalence.
\end{itemize}
We begin by relating instantons with holomorphic bundles.
For details see \cite{Buc93}.
The space $\overline{X_q}\setminus\{x_\infty\}$ is isomorphic
to the blow up $\Cc_q$
of $\Cc^2$ at $q$ points $x_1,\ldots,x_q$,
and given a self-dual
connection $A$ on $E\to X_q$,
its $(0,1)$ part $A^{0,1}$ is a holomorphic structure
on the bundle $E\to\Cc_q$. Denote by $\Ee_A$ the resulting
holomorphic bundle over $\Cc_q$. 
$\Cc_q$ can be compactified by adding a line $L_\infty$
at infinity.  The resulting compact surface is the
blow-up $\CP_q^2$ of $\CP^2$ at $q$ points.
Then $\Ee_A$ extends to a 
holomorphic bundle over $\CP^2_q$, and
the framing of the connection $A$ at $x_\infty\in X_q$
induces a framing of $\Ee$ at $L_\infty$.
In \cite{Buc93}, \cite{Mat00}, it was
shown that this correspondence is a real analytic isomorphism
betwen moduli spaces:

\begin{defi}
Let $p:\overline{\CP^2_q}\to X_q$ be the blow-up map
sending $L_\infty$ to $x_\infty$. Let $\mathcal M_k^r(\CP^2_q)$
be the quotient by ${\rm Aut}(p^*E)$ of the space of
pairs $(\alpha^{0,1},\phi)$, where $\alpha^{0,1}$ is a holomorphic
structure on $p^*E\to\CP_q^2$, inducing a holomorphic bundle
$\Ee$ trivial at $L_\infty$,
and $\phi:\Ee|_{L_\infty}\to\Oo_{L_\infty}^{\oplus r}$ 
is a trivialization. 
\end{defi}

\begin{teor}
Fix an isomorphism $h:E_{x_\infty}\to\Cc^r$ and define
a map $\Psi:\Mm_k^r(X_q)\to\mathcal M_k^r(\CP^2_q)$ by
$\Psi([A])=[p^*A^{0,1},p^*h]$. Then
$\Psi$ is an isomorphism of real analytic moduli spaces.
\end{teor}

The isomorphism $\Psi$ has the following properties:

\begin{prop}\label{prop-app-oplus}
Let $A_J=\bigoplus_a\lambda_a$ be the connection associated
with the $q$-tuple $J$. For each $a$, $\lambda_a$ induces
a holomorphic structure on $L_a$, defining a holomorphic
line bundle $\mathcal L_a$. Then
the holomorphic bundle $\Ee_{A_J}$ associated with the
connection $A_J$ is isomorphic to the 
direct sum $\bigoplus_a\mathcal L_a$. 
\end{prop}
\begin{proof}
It follows from $p^*A^{0,1}=\bigoplus_ap^*\lambda_a^{0,1}$
\end{proof}

\begin{prop}\label{prop-app-g}
Let $g\in \Gg/\Gg_0\cong{\rm Aut}\,E_{x_\infty}\cong SU(r)$.
Then $\Psi([g^{-1}Ag])=[p^*A^{0,1},(p^*g)\circ(p^*h)]$
($p^*g\in {\rm Aut}\,( \Oo_{L_\infty}^{\oplus r})$).
\end{prop}
\begin{proof}
It follows from $p^*(g^{-1}A^{0,1}g)=(p^*g)^{-1}(p^*A^{0,1})p^*g$.
\end{proof}

Let $L_i\subset\CP_q^2$, $i=1,\ldots,q$, denote the exceptional divisors.
Then (see \cite{San03})

\begin{teor}\label{app-t4.4}
Let $I\subset\{1,\ldots,q\}$ and write $|I|=\#I$.
Let $\pi_I:\CP^2_q\to \CP^2_{|I|}$ be the blow up at points
$x_j$, $j\notin I$. $\pi_I$ induces embeddings
\begin{equation}\label{eq-app-pi}
\pi_I^*:\mathcal M_k^r(\CP^2_{|I|})\to\mathcal M_k^r(\CP^2_q)  
\end{equation}
defined by $\pi_I^*([\alpha^{0,1},\phi])=[\pi_I^*\alpha^{0,1},\phi]$,
whose images are the bundles
\begin{equation}\label{eq-app-impi}
\pi_I^*(\mathcal M_k^r(\CP^2_{|I|}))=
\left\{[\Ee,\phi]\,|\,\Ee|_{L_j}\mbox{ is trivial for }j\notin I\right\}
\end{equation}
\end{teor}

The inverse map 
$(\pi_I^*)^{-1}:\pi_I^*\mathcal M_k^r(\CP^2_{|I|})\to\mathcal M_k^r(\CP^2_{|I|})$
can be described as follows: 
$\displaystyle\CP_q^2\setminus\bigcup_{j\notin I}L_j$ is biholomorphic
to $\CP^2_{|I|}$ minus $q-|I|$ points. Given a bundle over $\CP_q^2$
restrict it to $\CP_q^2\setminus\bigcup L_j$. Then this bundle extends
in a unique way to give a bundle over $\CP^2_{|I|}$. A corollary of
this description is

\begin{prop}\label{prop-app-Or}
To a $q$-tuple $J=(j_1,\ldots,j_q)$ such that $j_i=0$
whenever $i\notin I$ we can associate a $(q-|I|)$-tuple
$\tilde J$ formed by the entries $j_i$ of $J$ with $i\in I$.
Then $\pi_I^*Or(\tilde J)=Or(J)$. In particular
$Or(J)\subset\pi_I^*(\Mm_k^r(\CP^2_{|I|}))$.
\end{prop}
\begin{proof}
From proposition \ref{prop-app-oplus}
it follows that $\pi_I^*Or(\tilde J)\subset Or(J)$.
Now, from equation \ref{eq-app-impi}, 
$Or(J)\subset\pi_I^*(\Mm_k^r(\CP^2_{|I|}))$. Then
$(\pi_I^*)^{-1} Or(J)\subset Or(\tilde J)$ which completes
the proof.
\end{proof}

A direct consequence of proposition \ref{prop-app-g} is

\begin{prop}\label{prop-app-comm}
Let $A\in\Mm_k(X_{|I|})$. Then $\pi_I^*(Or(A))=Or(\pi_I^*A)$
\end{prop}

The last result of this section is

\begin{teor}\label{teor1234}
Let $f_c:\Mm_k(X_q)\to BU(k)$ be the classifying map of
$[{\rm Ker}\,\Dd_c]$.
Then the composition
\[
\xymatrix{
BU(k)\cong\Mm_k(S^4)\ar[r]^-{\pi_\emptyset^*}&
\Mm_k(X_q)\ar[r]^-{f_c}&BU(k)
}\]
is a homotopy equivalence.
\end{teor}
\begin{proof}
We first prove the result for $k=1$. Let 
$\lambda$ be a self-dual (irreducible) connection
on a rank 2 charge 1 bundle $E_2$ over $S^4$. Then we 
can define a self-dual connection
$A=\left(\begin{smallmatrix}\lambda&0\\0&0
\end{smallmatrix}\right)$
on $E=E_2\oplus\Oo^{r-2}$. 
The stabilizer
of $A$ is the subgroup $U(1)\times U(r-2)\subset U(r)$
of matrices
\[
g=\begin{pmatrix}e^{i\theta}\\&e^{i\theta}\\&&g_{22}
\end{pmatrix}
\]
where $g_{22}\in U(r-2)$. From proposition
\ref{prop-app-comm} $Or(\pi_\emptyset^*A)=
\pi_\emptyset^*(Or(A))$.
This shows that the maps
$f_{c}:Or(\pi_\emptyset^*A)\to BU(1)$ are homotopy equivalences.
The result of the theorem for $k=1$ then follows since
the inclusion $Or(A)\to\Mm_1(S^4)$ is a homotopy
equivalence.

Now we prove the general case. 
Direct sum induces a map
$\Mm_{k_1}^{r_1}(X_q)\times\Mm_{k_2}^{r_2}(X_q)
\to\Mm_{k_1+k_2}^{r_1+r_2}(X_q)$
which is well defined when $r_1,r_2\to\infty$
(see theorem \ref{teor-app-h} in the appendices).
We have the commutative diagram
\[
\xymatrix{
\Mm_1(X_0)\times\ldots\times\Mm_1(X_0)\ar[r]^-{h}
\ar[d]^{\pi_\emptyset^*\times\ldots\times\pi_\emptyset^*}&
\Mm_k(X_0)\ar[d]^{\pi_\emptyset^*}\\
\Mm_1(X_q)\times\ldots\times\Mm_1(X_q)\ar[r]^-{h}
\ar[d]^{f_c\times\ldots\times f_c}&
\Mm_k(X_q)\ar[d]^{f_c}\\
BU(1)\times\ldots\times BU(1)\ar[r]^-{\oplus}&
BU(k)}
\]
We showed already that the left vertical maps give a homotopy
equivalence. The bottom map (Whitney sum) is surjective in
homology. It follows that the composition
$f_c\circ\pi_\emptyset^*$ of the right vertical maps is surjective in
homology. Hence it must be an isomorphism in homology.
Since $BU(k)$ is simply connected, $f_c\circ\pi_\emptyset^*$
is a homotopy equivalence.
\end{proof}

\section{Applications}

The objective of this section is to prove theorems \ref{t3} and \ref{t2}.
They will be consequences of the theorem

\begin{teor}\label{teor41}
Let $\imath:\Sigma_{q,k}\to\Mm_k(X_q)$ be the inclusion and let
\[
K={\rm Ker}\left(\,\imath^*\circ f^*:H^*\left(
{\textstyle\prod} BU(k)\right)\to H^*(\Sigma_{q,k})\,\right)
\]
Then, for $k=1,2$
the map $\imath^*:H^*(\Mm_k(X_q))\to H^*(\Sigma_{q,k})$ is injective and the map
\[
f^*:H^*\left({\textstyle\prod} BU(k)\right)/K\to H^*(\Mm_k(X_q))
\]
is an isomorphism.
\end{teor}

We begin by computing $K$. Recall that $\Sigma_{q,k}=\coprod Or(J)$
and we have $Or(J)=\prod_{i=1}^q\left(\,BU(j_i)\times BU(j_i)\,\right)$.
We introduce some notation:
\begin{align*}
&H^*\left(\coprod_{|J|=k}Or(J)\right)\cong
\bigoplus_{|J|=k}
\left[\J c_n^{i+},
{\J c_n^{i-}};\,1\leq i\leq q,\,1\leq n\leq j_l\right]\\
&H^*\left(\prod_{l=0}^qBU(k)\right)\cong
\Zz\left[C_n^l;\,0\leq l\leq q,n\geq1\right]
\end{align*}

It will be convenient to make a change of variable. With power series 
notation ($C^l=1+C_1^l+\ldots$, ${\J c^{i\pm}}=1+{\J c^{i\pm}_1}+\ldots$) let
\begin{align*}
&C=C^0&&{\J c^i}={\J c^{i+}}\\
&S^i=(C^0)^{-1}C^i&&{\J s^i}=({\J c^{i+}})^{-1}{\J c^{i-}}
\end{align*}

Then
\begin{align}\label{S->s}
&f^*(C)=\sum_{|J|=k}{\J c^{1}}\ldots{\J c^{q}}&
&f^*(S^i)=\sum_{|J|=k}{\J s^i}
\end{align}

We will also use in this section the following notation:
Let $J_i$ be the $q$-tuple such that $j_l=0$ for $l\neq i$ and $j_i=1$.
Then let $J_{ij}=J_i+J_j$.

\begin{prop}
For $k=1$, $K$ is generated by $S_1^iS_1^j$, $i\neq j$. For $k=2$,
$K$ is generated by
\begin{align*}
&C_1S_1^iS_1^j+S_1^iS_2^j+S_2^iS_1^j\,,\ i\neq j\\
&C_2S_1^iS_1^j-S_2^iS_2^j\,,\ i\neq j\\
&S_1^iS_1^jS_1^k\,,\ i\neq j\neq k\\
&S_1^iS_1^jS_2^k\,,\ i\neq j\neq k
\end{align*}
\end{prop}
\begin{proof}
These elements are clearly in the kernel so we only have to
show the map is injective if we mod out by this ideal. For $k=1$ this
is clearly the case so let $k=2$.
Then, for each $i$, the composition
\[
\Zz[C_1,C_2,S_1^i,S_2^i]\to H^*\left({\textstyle\prod} BU(2)\right)\to 
H^*(\Sigma)\to Or(J_{ii})
\]
is an isomorphism. Hence the kernel is contained in the ideal $M_2$
generated by products
$S_{n_1}^{i_1}S_{n_2}^{i_2}$, $n_j=1,2$, $i_1\neq i_2$.
Now observe that
\[
S_2^i(C_1S_1^iS_1^j+S_1^iS_2^j+S_2^iS_1^j)+
S_1^i(C_2S_1^iS_1^j-S_2^iS_2^j)=
S_1^j((S_2^i)^2+C_1S_1^iS_2^i+C_2(S_1^i)^2)
\]
It follows that any element $x\in M_2$ can be written as
\[
x=\sum_{i\neq j}(P_{ij}S_1^iS_1^j+Q_{ij}S_1^iS_2^j)+k
\]
where $P_{ij},Q_{ij}\in\Zz[C_1,C_2,S_1^i,S_1^j]$ and $k\in K$.
We want to show that if $f^*(x)=0$ then
$P_{ij}=Q_{ij}=0$. Fixing $i,j$, $i\neq j$, 
and letting $f_{ij}^*$
be the composition $H^*\left(\prod BU(2)\right)\to
H^*(\Sigma)\to Or(J_{ij})$ we have
\[
f^*_{ij}(x)=(f^*_{ij}(P_{ij})+f^*_{ij}(Q_{ij})c_1^j)s_1^is_1^j=0
\]
But both $f^*_{ij}(P_{ij})$ and $f^*_{ij}(Q_{ij})$
are symetric in $c_1^i,c_1^j$ hence
we must have $f^*_{ij}(P_{ij})=f^*_{ij}(Q_{ij})=0$.
This clearly implies $P_{ij}=Q_{ij}=0$.
\end{proof}

Before we prove theorem \ref{teor41} we need the following results
proven in \cite{San03}:

\begin{teor}\label{t4.3}
Consider the ideals $K_i,K_{ij}$ generated by
\begin{align*}
&K_i=\left\langle k_1,k_2\right\rangle\subset\Zz[a_1,a_2,k_1,k_2]\cong 
H^*(BU(2)^{\times 2})\\
&K_{ij}=\left\langle x_1x_2\right\rangle\subset\Zz[x_1,x_2,x_3,x_4]
\cong H^*(BU(1)^{\times 4})
\end{align*}
Then, we have the isomorphism of abelian groups
\[
H^*(\Mm_2(X_q))\cong\Zz[a_1,a_2]\oplus\bigoplus_iK_i\oplus\bigoplus_{i<j}K_{ij}
\]
\end{teor}

The next theorem uses the notation in theorem \ref{app-t4.4}
with $I=\emptyset,\{i\},\{i,j\}$:

\begin{teor}\label{t4.4}
Let $C\subset\Mm_2(X_2)\cong\mathcal M_2(\CP^2_2)$ be given by
\begin{equation}\label{eqc}
C=\left\{[\Ee,\phi]\,|\,\Ee\mbox{ is trivial in }\CP^2_2\setminus
(L_1\cup L_2)\right\}
\end{equation}
There is a homotopy equivalence $g:C\to\Mm_1(X_1)\times\Mm_1(X_1)$ and
\begin{enumerate}
\item Let $\imath_\emptyset^*:\pi_\emptyset^*\Mm_2(X_0)\to\Mm_2(X_q)$ be the inclusion.
Then $\imath_\emptyset^*:\Zz[a_1,a_2]\to H^*(\pi_\emptyset^*\Mm_2)$ is
an isomorphism.
\item Let $\imath_i:\pi_i^*\Mm_2(X_1)\to\Mm_2(X_q)$ be the inclusion.
Then $\imath^*_i:\Zz[a_1,a_2]\oplus K_i\to H^*(\pi_i^*\Mm_2)$
is an isomorphism and
$\imath^*_i(K_{kl})=\imath^*_i(K_j)=0$ for any $k,l$ and any $j\neq i$.
\item Let $\imath_{ij}:\pi_{ij}^*C\to\Mm_2(X_q)$ be the inclusion.
Then the map $\imath^*_{ij}:K_{ij}\to H^*(\pi_{ij}^*C)$ is injective
and $\imath^*_{ij}(K_{kl})=0$ for different sets
$\{k,l\}\neq\{i,j\}$.
\end{enumerate}
\end{teor}

\begin{prop}
$Or(J_{ii})\subset\pi_i^*\Mm_2$, $Or(J_{ij})\subset\pi_{ij}^*C$ ($i\neq j$), and these inclusions
are homotopy equivalences.
\end{prop}
\begin{proof}
From proposition \ref{prop-app-Or} we have 
$\pi_i^*(\Sigma_1)=Or(J_{ii})$ and from theorem \ref{teorA2},
it follows that the inclusion $Or(J_{ii})\to\pi_i^*\Mm_2$
is a homotopy equivalence.
Since $Or(1,1)\subset C$ ($Or(1,1)$ is the orbit associated with
$J=(1,1)$), proposition \ref{prop-app-Or} implies
$Or(J_{ij})=\pi_{ij}^*Or(1,1)\subset\pi_{ij}^*C$.
To show that this inclusion
is a homotopy equivalence we need to look at the map
$g:C\to\Mm_1(X_1)\times\Mm_1(X_1)$. This map is defined as follows
(compare with the discussion after theorem \ref{app-t4.4}): 
the spaces  $\CP^2_2\setminus L_i$ ($i=1,2$)
are biholomorphic to
$\CP_1^2$ minus one point.
The bundles
$\Ee|_{\CP^2_2\setminus L_i}$ extend uniquely
to give bundles
$\Ee_i$ over $\CP^2_1$ $(i=1,2)$. Then $g$ is 
defined by $g(\Ee)=(\Ee_1,\Ee_2)$.
We have a commutative diagram
\[
\xymatrix{
\pi_{ij}^*\Mm_2(X_2)&
\pi_{ij}^*C\ar[l]&
C\ar[l]^-{\cong}_-{\pi_{ij}^*}\ar[r]^-{g}_-{\cong}&
\Mm_1(X_1)\times\Mm_1(X_1)\\
&
Or(J_{ij})\ar[u]\ar[lu]&
Or(1,1)\ar[u]\ar[l]^-{\cong}_-{\pi_{ij}^*}\ar[r]^-{g}&
\Sigma_{1,1}\times\Sigma_{1,1}\ar[u]^{\cong}
}\]
We claim that $g:Or(1,1)\to\Sigma_{1,1}\times\Sigma_{1,1}$
is a homotopy equivalence. 
It will follow that
the inclusion $Or(J_{ij})\to\pi_{ij}^*C$ is a homotopy equivalence.
To prove the claim we recall that
\begin{align*}
&Or(1,1)     \cong\frac{SU(r)}{U(1)^{\times 4}\times U(r-4)}&
&\Sigma_{1,1}\cong\frac{SU(r)}{U(1)^{\times 2}\times U(r-2)}
\end{align*}
Then it follows from proposition \ref{prop-app-g} that the map
$Or(1,1)\to\Sigma_{1,1}$ is the natural projection map. When
$r\to\infty$ we see that $g:Or(1,1)\to\Sigma_{1,1}\times\Sigma_{1,1}$
is a homotopy equivalence. This concludes the proof.
\end{proof}

Now we prove theorem \ref{teor41}:

\begin{proof}
First we observe that the map 
$\imath^*:H^*(\Mm_2(X_q))\to H^*(\Sigma)$ is injective.
This follows easily from theorems \ref{t4.3} and \ref{t4.4}.
Hence we only need to show that $f^*:H^*(\prod BU(2))\to
H^*(\Mm_2)$ is surjective.
We will divide the proof into several steps.
\begin{enumerate}
\item
Let $I_{ij}\subset\Zz[C_1,C_2,S_1^l,S_2^l]$ be the ideal
generated by $S_1^iS_1^j$ and $S_1^iS_2^j$.
We claim that $f^*(I_{ij})\subset K_{ij}$. Let
$y\in I_{ij}$. Then write
$f^*(y)=x_\emptyset+\sum x_i+\sum x_{ij}$ with
$x_\emptyset\in\Zz[a_1,a_2]$, $x_i\in K_i$ 
and $x_{ij}\in K_{ij}$.
Then, for any $l$,
$\imath_l^*f^*(y)=0=\imath_l^*(x_\emptyset+x_l)$ 
so $x_\emptyset=x_l=0$.
Then, for any pair $\{k,l\}\neq\{i,j\}$ we have
$\imath_{k,l}^*f^*(y)=0=\imath_{k,l}^*(x+x_k+x_l+x_{kl})=
\imath_{k,l}^*(x_{kl})$ hence $x_{kl}=0$. 
This completes the proof.
\item 
We claim that $f^*(I_{ij})=K_{ij}$. First we
observe that the image of the map
$I_{ij}\to\Zz[c_1^i,c_1^j,s_1^i,s_1^j]$ is precisely the ideal
generated by $s_1^is_1^j$. It is clearly contained in that ideal so
we only have to show surjectivity. Recall (equation \ref{S->s}) that
$S^i_n\mapsto s^i_n$. Also, since $c_n^i=0$ for $n>1$ it follows that
$s_n^i=(-c_1^i)^{n-1}s_1^i$. This shows surjectivity.
Now, since
$f^*(I_{ij})\cong K_{ij}$ and $s_1^i,s_1^j$ are
both primitive elements in $H^*(Or(J_{ij}))$ it follows that
$f^*(I_{ij})=K_{ij}$.
\item
Now let $I_i\subset\Zz[C_1,C_2,S_1^i,S_2^i]$ be the ideal
generated by $S_1^i,S_2^i$.
We claim that $\imath_i^*f^*(I_i)=\imath_i^*(K_i)$. Let $y\in I_i$.
Then write $f^*(y)=x_\emptyset+\sum x_i+\sum x_{ij}$. 
Then, for any $l\neq i$,
$\imath_l^*f^*(y)=0=\imath_l^*(x_\emptyset+x_l)$ so 
$x_\emptyset=x_l=0$. Hence
$\imath^*f^*(I_i)\subset\imath^*K_i$. Now, since
$\imath^*f^*(I_i)\cong \imath_i^*(K_i)$ and both
$S_1^i,S_2^i$ are sent to primitive elements, it follows that
$\imath_i^*f^*(I_i)=\imath_i^*(K_i)$.
\item
Let $x\in K_i$. We claim that there is $y\in H^*(\prod BU(2))$
such that $f^*(y)=x$. Proof:
choose $y_i\in I_i$ such that $\imath^*_if^*(y_i)=\imath_i^*(x)$. 
Then $f^*(y_i)=x_i+\sum x_{kl}$. So we only have to 
choose $y_{kl}$ such that $f^*(y_{kl})=x_{kl}$ and define
$y=y_i-\sum x_{kl}$.
\item Let $x\in\Zz[a_1,a_2]$. We claim that
there is $y\in H^*(\prod BU(2))$
such that $f^*(y)=x$. From theorem \ref{teor1234}, the map
$\imath_\emptyset^*:\Zz[C_1,C_2]\to H^*(S^4)$
is an isomorphism hence we can find 
$y_\emptyset\in\Zz[C_1,C_2]$ such that
$\imath_\emptyset^*f^*(y_\emptyset)=\imath_\emptyset^*(x)$.
Then $f^*(y_\emptyset)=x+\sum_ix_i+\sum x_{ij}$
so we just choose $y_i\mapsto x_i$, $y_{ij}\mapsto x_{ij}$
and let $y=y_\emptyset-\sum y_i-\sum y_{ij}$.
\end{enumerate}
\end{proof}

\appendix

\section{Monads}

Here we show that the inclusion $\Sigma_{1,k}^\infty\to\Mm_k^\infty(X_1)$ 
is a homotopy equivalence.
We will use the monad description of self-dual connections on $X_1$
(see \cite{Kin89}): let $U,W$ be $k$ dimensional vector spaces
and consider the space $\bar R$ of
configurations $(a_1,a_2,d,b,c)$ where
$a_1,a_2\in{\rm Hom}(U,W)$, $d\in{\rm Hom}(W,U)$, 
$b\in{\rm Hom}(\Cc^r,W)$ and $c\in{\rm Hom}(U,\Cc^r)$,
satisfying the integrability condition
$a_1da_2-a_2da_1+bc=0$.
The group ${\rm Aut}(W)\times{\rm Aut}(U)$ acts in $\bar R$ by
\[
(g,h)\cdot(a_1,a_2,d,b,c)=(ga_1h^{-1},ga_2h^{-1},hdg^{-1},gb,ch^{-1})
\]
Then there is an open subset of $\bar R$, $R$
(the non-degenerate configurations)
such that $\Mm_k^r(X_1)$ is isomorphic to the quotient of $R$
by the action of ${\rm Aut}(W)\times{\rm Aut}(U)$.

\begin{prop}
Let $M_0$ be the set of equivalence classes of configurations of the
form $(0,0,0,b,c)$. Then, as a subspace of $\Mm_k^r(X_1)\cong\mathcal M_k^r(\CP^2_1)$, 
$M_0=\Sigma_{1,k}^r$.
\end{prop}
\begin{proof}
A self dual connection $A\in M_0\subset\Mm_k^r(X_1)$ induces a
holomorphic bundle $\mathcal E$
on $\CP_1^2$ given by the homology of the 
sequence ($\mathcal E={\rm Ker}\,B/{\rm Im}\,A$)
\begin{multline}\label{seqa}
\xymatrix{
U\otimes\Oo(-L_\infty)\oplus W\otimes\Oo(-L_\infty+L)\ar[r]^-{A}&
V\otimes\Oo\ar[r]^-{B}&}\\
\xymatrix{\ar[r]^-{B}&
W\otimes\Oo(L_\infty)\oplus U\otimes\Oo(L_\infty-L)
}\end{multline}
where
\[
A=\begin{pmatrix}0&s_1\Id\\x_1\Id&0\\0&s_2\Id\\x_2\Id&0\\cx_3&0\end{pmatrix}\,,\ 
B=\begin{pmatrix}x_2\Id&0&-x_1\Id&0&bx_3\\0&s_2\Id&0&-s_1\Id&0\end{pmatrix}
\]
Here $\Id$ is the $k\times k$ identity matrix,
$x_1,x_2,x_3$ are sections spanning $H^0(\Oo(L_\infty))$ 
with $x_3=0$ in $L_\infty$ and
$s_1,s_2$ are sections spanning
$H^0(\Oo(L_\infty-L))$ such that $s_1x_2=s_2x_1$. 
Let $\Ee_{1i},\Ee_{2i}\subset{\rm Ker}\,B$ be the rank 2 bundles
\begin{align*}
&\Ee_{1i}=\left\{
\begin{pmatrix}0&\alpha s_1\Id&0&\alpha s_2\Id&\beta c\end{pmatrix}_i
\,;\,\alpha,\beta\in\Cc\right\}\\
&\Ee_{2i}=\left\{
\begin{pmatrix}v_2\Id&0&-v_1\Id&0&v_3\frac{\bar b^t}{|b|^2}\end{pmatrix}_i
\,;\,v_1x_1+v_2x_2+v_3x_3=0\right\}
\end{align*}
(the subscript $i$ denotes the $i$th line of the matrix)
and let $\Ee_{3}=\{(0,0,0,0,v)|bv=\bar c^tv=0\}$.
Then ${\rm Ker}\,B=\bigoplus_i(\Ee_{1i}\oplus\Ee_{2i})\oplus\Ee_{3}$.
Sequence \ref{seqa} splits into $k$ sequences whose homology is the homology of the sequences
\[
\xymatrix@R=0em@C=6em{
\Oo(-L_\infty)\ar[r]^-{(0,x_1,0,x_2,c_ix_3)}&\Ee_{1i}\ar[r]&0\\
\Oo(-L_\infty+L)\ar[r]^-{(s_1,0,s_2,0,0)}&\Ee_{2i}\ar[r]&0
}\]
This clearly shows that the homology of sequence \ref{seqa} is a sum of line
bundles $\Ee\cong\bigoplus_i(\Oo(L)\oplus\Oo(-L))\oplus\Oo^{r-2k}$.
Hence $M_0\subset\Sigma_{1,k}^r$. To show the opposite direction we
observe that $M_0$ is the subspace of pairs
$(V_1,V_2)\in Grass(k,\Cc^r)\times Grass(k,\Cc^r)$ with $V_1\perp V_2$
hence it has dimension $4rk-6k^2$. Hence, $M_0,\Sigma_{1,k}^r$ are 
compact manifolds of the same dimension. It follows that $\Sigma_{1,k}^r=M_0$. 
\end{proof}

\begin{teor}\label{teorA2}
When $r\to\infty$,
the inclusion $\imath:\Sigma_{1,k}^\infty\to\Mm_k^\infty(X_1)$ is a homotopy equivalence.
\end{teor}
\begin{proof}
Consider the principal $Aut(W)\times Aut(V)$ bundle $R\to\Mm_k^\infty(X_1)$.
Consider the pullback of this bundle to $M_0$. We have the diagram
\[
\xymatrix{
\imath^*R\ar[r]\ar[d]&R\ar[d]\\
M_0\ar[r]&\Mm_k^\infty
}\]
In \cite{BrSa97} it was shown that $R$ becomes contractible 
when $r\to\infty$ and
the same argument shows that  $\imath^*R$ (the set of configurations
of the form $(0,0,0,b,c)$) is also contractible.
Applying the five lemma to the long exact sequences 
of homotopy groups associated
to these principal bundles it follows that the map $\imath_*$
is an isomorphism in all homotopy groups, hence $\imath$ is a homotopy
equivalence.
\end{proof}

\section{Direct sum}

Here we will define a map
$\Mm_{k_1}^{r_1}(X_q)\times\Mm_{k_2}^{r_2}(X_q)
\to\Mm_{k_1+k_2}^{r_1+r_2}(X_q)$ induced
by direct sum and show that this map is well
defined when $r_1,r_2\to\infty$.

For each pair $(r,k)$ fix an $SU(r)$ bundle
$E_{r,k}$ over $X_q$ with
$c_2(E_{r,k})=k$. Then, to each isomorphism
$\phi:E_{k_1,r_1}\oplus E_{k_2,r_2}\to
E_{k_1+k_2,r_1+r_2}$ we can associate a map
\[
h_\phi:\Mm_{k_1}^{r_1}(X_q)\times\Mm_{k_2}^{r_2}(X_q)
\to\Mm_{k_1+k_2}^{r_1+r_2}(X_q)
\]
that sends a pair of connections $(A_1,A_2)$
to $\phi^*(A_1\oplus A_2)$. Since we act
on the connections with the gauge group $\Gg_0$,
$h_\phi$ only depends on the value of $\phi$ at
the base point $x_\infty\in X_q$.
The following result is then an easy consequence
of the connectedness of $SU(r_1+r_2)$:

\begin{prop}\label{qwerty}
Any two maps $h_{\phi_1},h_{\phi_2}$ are homotopic.
\end{prop}

Sum with a trivial rank $R$ bundle induces a map
$\Mm_k^r\to\Mm_k^{r+R}$. Then

\begin{teor}\label{teor-app-h}
The following diagram is homotopy commutative
\[
\xymatrix{
\Mm_{k_1}^{r_1}\times\Mm_{k_2}^{r_2}\ar[r]^-{h}\ar[d]&
\Mm_{k_1+k_2}^{r_1+r_2}\ar[d]\\
\Mm_{k_1}^{r_1+R}\times\Mm_{k_2}^{r_2+R}\ar[r]^-{h}&
\Mm_{k_1+k_2}^{r_1+r_2+2R}
}\]
Hence the map $h:\Mm_{k_1}^\infty\times\Mm_{k_2}^\infty
\to\Mm_{k_1+k_2}^\infty$ is well defined.
\end{teor}
\begin{proof}
It follows directly from proposition \ref{qwerty}.
\end{proof}

\bibliography{bi}

\end{document}